\documentclass[12pt,twoside]{article}
\usepackage{amsmath}
\usepackage{mathrsfs}
\usepackage{amsfonts}
\usepackage{amssymb}

\usepackage{amsfonts,amsmath,amsthm}
\allowdisplaybreaks

\numberwithin{equation}{section} \textwidth 15.5 true cm
\begin{document}
  \title{{\bf  LPS's Criterion for Incompressible Nematic Liquid
  Crystal Flows
  }}
\author{\it\small
Qing Chen\\\it\small Department of Mathematics and
Physics\\\it\small
Xiamen University of Technology, Xiamen 361024, China\\[3mm]\it\small
Zhong Tan\\\it\small School of Mathematical Sciences, Xiamen
University,
 Fujian 361005, China\\\it\small
Guochun Wu\thanks{Corresponding author. Email address:
guochunwu@126com.}\\\it\small School of Mathematical Sciences,
Xiamen University,
 Fujian 361005, China}

  \date{}
  \maketitle
  \begin{abstract}{\small
    In this paper we derive LPS's criterion for the breakdown of classical solutions
    to the incompressible nematic liquid crystal flow, a simplified version of Ericksen-Leslie
    system modeling the hydrodynamic evolution of nematic liquid crystals in $\mathbb
    R^3$. We show that if $0<T<+\infty$} is the maximal time
    interval for the unique smooth solution $u\in
    C^\infty([0,T),\mathbb
    R^3)$, then $|u|+|\nabla d|\notin L^q([0,T],L^p(\mathbb
    R^3))$, where $p$  and $q$ safisfy the
    Ladyzhenskaya-Prodi-Serrin's
    condition: $\frac{3}{p}+\frac{2}{q}=1$ and $p\in(3,+\infty]$.
     \\
 \\
{\small {\bf Keywords} Incompressible nematic liquid crystal flow;
Ladyzhenskaya-Prodi-Serrin's criterion.}\smallskip
\\
{\small {\bf MSC2010} 35Q35; 76D03.}

  \end{abstract}

\section{Introduction}

$\ \ \ \ $We consider the following hydrodynamic system modeling the
flow of liquid crystal materials in dimension three (see [2,3,10,12]
and references therein):
$$ u_t+u\cdot\nabla u-\nu\triangle u+\nabla P=-\triangle
d\cdot\nabla
d,\eqno(1.1a)$$ $$\partial_t d+u\cdot\nabla d=\triangle d +|\nabla
d|^2d,\eqno (1.1b)
$$$$\nabla\cdot u=0\ \ |d|=1,\eqno(1.1c)$$
for $(t,x)\in[0,+\infty)\times\mathbb{R}^{3}$. Here $u:\mathbb
R^3\rightarrow \mathbb R^3$ represents the velocity field of the
incompressible viscous fluid, $\nu>0$ is the kinematic viscosity,
$P:\mathbb R^3\rightarrow \mathbb R$ represents the pressure
function, and $d:\mathbb R^3\rightarrow \mathbb S^2$ represents the
macroscopic average of the nematic liquid crystal orientation field.
We are interested in the Cauchy problem (1.1) with the initial value
$$(u(0,x),d(0,x))=(u_0(x),d_0(x))\eqno (1.2)$$
satisfying the following compatibility condition:
$$\nabla \cdot u_0(x)=0,\ \ |d_0(x)|=1,\ \ \lim\limits_{|x|\rightarrow\infty} d_0(x)=a\in\mathbb S^2,\eqno (1.3)$$
where $a$ is a given unit vector.

The above system is a simplified version of the Ericksen-Leslie
model, which reduces to the Ossen-Frank model in the static case,
for the hydrodynamics of nematic liquid crystals developed during
the period of 1958 through 1968 [2,3,10]. It is a macroscopic
continuum description of the time evolution of the materials under
the influence of both the flow field $u(x,t)$, the macroscopic
description of the microscopic orientation configurations $d(x,t)$
of rod-like liquid crystals. Roughly speaking, the system (1.1) is a
coupling between the non-homogeneous Navier-Stokes equation and the
transported flow harmonic maps. Due to the physical importance and
mathematical challenges, the study on nematic liquid crystals has
attracted many physicists and mathematicians. The mathematical
analysis of the liquid crystal flows was initiated by Lin [11], Lin
and Liu in [12,13]. For any bounded smooth domain in $\mathbb R^2$,
Lin , Lin and Wang [14] have proved the global existence of
Leray-Hopf type weak solutions to system (1.1) which are smooth
everywhere except on finitely many time slices (see [5] for the
whole space). The uniqueness of weak solutions in two dimension was
studied by [15,20]. Recently, Hong and Xin [6] studied the global
existence for general Ericksen-Leslie system in dimension two.
However, the global existence of weak solutions to the
incompressible nematic liquid crystal flow equation (1.1) in three
dimension with large initial data is still an outstanding open
question.

In this paper, we are interested in an optimal characterization on
the maximal interval $T$ that is scaling invariant. So let us first
introduce the following definition:\\
\\{\bf Definition 1.1.} For $1\leq p,q\leq\infty$, we say a function
$f=f(t,x):[0,T]\times \mathbb R^3\rightarrow \mathbb R$ is in
$L^q([0,T],L^p(\mathbb R^3))$, if
$$\begin {array}{rl}\|f\|_{L^q([0,T],L^p(\mathbb R^3))}&=(\int_0^T\|f(t,\cdot)\|^q_{L^p(\mathbb R^3)}dt)^\frac{1}{q},\ \ 1\leq q<\infty
\\&=ess\sup\limits_{t\in[0,T]}\|f(t,\cdot)\|_{L^p(\mathbb R^3)},\ \ q=\infty\end {array}$$
is finite. If $p=q$, then we simply write $\|f\|_{L^p([0,T],\mathbb
R^3)}$ for $\|f\|_{L^p([0,T],L^p(\mathbb R^3))}$.\\

We will consider the short time classical solution to (1.1) and
address the Ladyzhenskaya-Prodi-Serrin's criterion that
characterizes the first time  finite singular time. The local
well-posedness of the Cauchy problem of system (1.1) is rather
standard (see [5,7,14]). More precisely, if the initial velocity
$u_0\in H^s(\mathbb R ^3,\mathbb R ^3)$ with $\nabla\cdot u_0=0$ and
$d_0-a\in H^{s+1}(\mathbb R^2,\mathbb S^2)$ such that system (1.1)
has a unique, classical solution $(u,d)$ in $[0,T_0)\times\mathbb
R^3$ satisfying $$\begin {array}{rl}&u\in C([0,T),H^s(\mathbb
R^3))\cap C^1([0,T),H^{s-2}(\mathbb R^3))\ \ and\\&d-a\in
C([0,T),H^{s+1}(\mathbb R^3,\mathbb S^2))\cap
C^1([0,T),H^{s-1}(\mathbb R^3,\mathbb S^2)),
\end {array}\eqno (1.4)$$for any $0<T<T_0$. At present, there is no global-in-time existence theory for
classical solutions to system (1.1). Thus if we assume $T_*>0$ is
the maximum value such that (1.4) holds with $T_0=T_*$, we would
like to characterize such a $T_*$. Motivated by the famous work [1],
Huang and Wang [7] have obtained a $BKM$ type blow-up criterion (see
also [16]). However, the techniques involved in this paper are much
different from [7], which we believe that the result may have its
own interest.

When $d$ is a constant vector field, the system (1.1) becomes an
incompressible Navier-Stokes equation. Recall that the scaling
invariant space $L^q([0,T],L^p(\mathbb R^3))$, with $(p,q)$
satisfying
$$\frac{3}{p}+\frac{2}{q}=1,\eqno (1.5)$$
which has played an important role in the regularity issue of
Navier-Stokes equation. Leray [9] first established the existence of
a global weak solution for Navier-Stokes equation, now called
Leray-Hopf weak solution, that satisfies an energy inequality:
$$\frac{1}{2}\|u(t)\|^2_{L^2(\mathbb R^3)}+\int^t_0\int_{\mathbb R^3}
|\nabla u(t,x)|^2dxds\leq \frac{1}{2}\|u_0\|^2_{L^2(\mathbb R^3)}.$$
Although the regularity issue for Leray-Hopf weak solutions of
Navier-Stokes equation remains open, it is well-known that both
uniqueness and smoothness for the class of weak solutions of
Navier-Stokes equation, in which $u\in L^q([0,T],L^p(\mathbb R^3))$
for some $p\in (3,+\infty]$ and $q\in[2,+\infty)$ satisfying
Ladyzhenskaya-Prodi-Serrin's condition (1.5) have been established
through works by Prodi [18], Serrin [19], and Ladyzhenskaya [8] in
1960s. On the other hand, for the end point case $p=3,q=+\infty$ ,
only until very recently Escauriaza et al. [4] have finally proved
the smoothness for weak solution $u\in L^\infty([0,T],L^3(\mathbb
R^3))$ of incompressible Navier-Stokes equation, $0<T\leq +\infty$.

Motivated by these results for the Navier-Stokes equation, we are
going to use scaling considerations for system (1.1) to guess which
spaces may be critical. We observe that system (1.1) is invariant by
the following transformation:
$$\hat u=lu(l^2t,lx),\ \ \hat P=l^2P(l^2t,lx),\ \ \hat d=d(l^2t,lx).$$
Thus $L^q([0,T],L^p(\mathbb R^3))$ is a critical space for
$(u,\nabla d)$ if $(p,q)$ satisfies the Ladyzhenskaya-Prodi-Serrin's
condition (1.5).

Our main
results are formulated as the following theorem:\\
\\{\bf Theorem 1.1.} For $u_0\in H^s(\mathbb R^3)$ with $\nabla\cdot u_0=0$ and $d_0-a\in H^{s+1}(\mathbb R^3)$
with $|d_0|=1$ for $s\geq 3$. Suppose that $(u,d)$ is a smooth
solution to the system (1.1)-(1.2), then for given $T>0$, $(u,d)$ is
smooth up to time $T$ provided that
$$\|u\|_{L^q([0,T],L^p(\mathbb R^3))}+\|\nabla d\|_{L^q([0,T],L^p(\mathbb R^3))}<+\infty,\eqno (1.7)$$
where $(p,q)$ satisfies the Ladyzhenskaya-Prodi-Serrin's condition
(1.5) and $p\in(3,+\infty]$.\\

{\bf Notations.} We denote by $L^p$, $W^{m,p}$ the usual Lebesgue
and Sobolev spaces on $\mathbb R^3$ and $H^m= W^{m,2}$, with norms
$\|\cdot\|_{L^p}$, $\|\cdot\|_{W^{m,p}}$ and $\|\cdot\|_{H^m}$
respectively. For the sake of conciseness, we do not distinguish
functional space when scalar-valued or vector-valued functions are
involved. We denote
$\nabla=\partial_x=(\partial_1,\partial_2,\partial_3)$, where
$\partial_i=\partial_{x_i}$, $\nabla_i=\partial_i$ and put
$\partial_x^lf=\nabla^lf=\nabla(\nabla^{l-1}f)$. We assume $C$ be a
positive generic constant throughout this paper that may vary at
different places and the integration domain $\mathbb R^3$ will be
always omitted without any ambiguity. Finally, $\langle
\cdot,\cdot\rangle $ denotes the inner-product in $L^2(\mathbb
R^3)$.\\
\\{\bf Remark 1.1.} It is standard that the condition (1.3) is
preserved by the flow. In fact, first notice that the divergence
free of the velocity field $u$ can be justified by the initial
assumption that $\nabla\cdot u_0=0$. Indeed, this can be easily and
formally observed by take $\nabla\cdot$ to the momentum equation.
Moreover, applying the maximum principle to the equation for
$|d|^2$, one also can easily see that $|d|=1$ under the initial
assumption that $|d_0|=1$.

\section{Proof of Theorem 1.1}

$\ \ \ \ $We prove our theorem in this section. Without loss of
generality, we assume $\nu=1$. The first bright idea to reduce many
complicated computations lies in that we just need to do the lowest
order and highest order energy estimates for the solutions. This is
motivated by the following observation:
$$\|f\|^2_{H^k}\leq C\|(f,\nabla^kf)\|^2_{L^2},\ \ \forall f\in H^k.\eqno (2.1)$$
This inequality (2.1) can be easily proved by combing Young's
inequality and Gagliardo-Nirenberg's inequality.
$$\|\nabla^if\|_{L^p}\leq C(p)\|f\|^\alpha_{L^q}\|\nabla^kf\|^{1-\alpha}_{L^r},\ \ \forall f\in H^k\eqno (2.2)$$
where
$\frac{1}{p}-\frac{i}{3}=\frac{1}{q}\alpha+(\frac{1}{r}-\frac{k}{3})(1-\alpha)$
with $i\leq k.$\\

Now we are in a position to prove our Theorem 1.1.\\
\\{\bf Proof of Theorem 1.1} First of all, we note that if $p=+\infty$, Theorem 1.1 has been proved in [16],
 thus let us concentrate on $p\in(3,+\infty)$.
Now for classical solutions to (1.1)-(1.2), one has the following
basic energy law:
$$\begin {array}{rl}\|u(t,\cdot)\|^2_{L^2}+\|\nabla d(t,\cdot)\|^2_{L^2}&+\int_0^t(\|\nabla u(s,\cdot)\|^2_{L^2}+\|\triangle d(s,\cdot)
+|\nabla d|^2d(s,\cdot)\|^2_{L^2})ds\\=&\|u_0\|^2_{L^2}+\|\nabla
d_0\|^2_{L^2},\ \ \forall\  t>0.\end {array}\eqno (2.3)$$

Let's concentrate on the case $s=3$. For each multi-index $\alpha$
with $|\alpha|\leq 3$, by applying $\partial^\alpha_x$ to (1.1a) and
$\partial^{\alpha+1}_x$ to (1.1b), multiplying them by
$\partial^\alpha_x u$, $\partial^{\alpha+1}_x d$ respectively and
then integrating them over $\mathbb R^3$, we have
$$\begin {array}{rl}&\frac{1}{2}\frac{d}{dt}\|\partial_x^\alpha (u,\nabla d)\|^2_{L^2}+\|\partial_x^\alpha (\nabla u,\triangle d)\|^2_{L^2}
\\=&-\langle \partial_x^\alpha (u\cdot\nabla u),\partial_x^\alpha u\rangle-\langle \partial_x^\alpha (\triangle d\cdot\nabla d),\partial_x^\alpha u\rangle
\\&-\langle \partial_x^{\alpha+1} (u\cdot\nabla d),\partial_x^{\alpha+1} d\rangle+\langle \partial_x^{\alpha+1} (|\nabla d|^2d),\partial_x^{\alpha+1} d\rangle
\\=&\sum\limits_{i=1}^4I_{|\alpha|,i}.\end {array}\eqno (2.4)$$
where $I_{|\alpha|,i}$ are the corresponding terms in the above
equation which will be estimated as follows. Now for $|\alpha|=1$ in
(2.4), integrating by parts and using the divergence free condition
$\nabla\cdot u=0$ and (2.2), we arrive at
$$\begin {array} {rl}|I_{1,1}|&=|\langle \partial_x^1 (u\cdot\nabla u),\partial_x^1 u\rangle|=|\langle \partial_x^1 u\cdot\nabla u,
\partial_x^1 u\rangle|\leq C\int_{\mathbb R^3} |\nabla u|^3dx\\&\leq C\|u\|_{L^p}^\frac{3p}{6+p}\|\nabla ^2u\|_{L^2}^\frac{18}{6+p}
\leq C\|u\|_{L^p}^\frac{2p}{p-3}\|u\|_{L^p}^\frac{p}{p-3}
+\frac{1}{16}\|\nabla ^2u\|^2_{L^2}.\end {array}\eqno (2.5)$$
Combining Cauchy's inequality, Sobolev's inequality and the fact
$|\nabla d|^2=-d\cdot\triangle d$ (since $|d|=1$) gives
$$\begin {array}{rl}|I_{1,2}|&=|\langle \partial_x^1 (\triangle d\cdot\nabla d),\partial_x^1 u\rangle|=
|\langle\triangle d\cdot\nabla d,\Delta u\rangle|\\&\leq
C\|\triangle d\|^3_{L^3}+\frac{1}{16}\|\nabla^2 u\|^2_{L^2}\\&\leq
C\|\nabla d\|_{L^p}^\frac{3p}{6+p}\|\nabla
^3d\|_{L^2}^\frac{18}{6+p}+\frac{1}{16}\|\nabla^2 u\|^2_{L^2}\\&
\leq C\|\nabla d\|_{L^p}^\frac{2p}{p-3}\|\nabla
d\|_{L^p}^\frac{p}{p-3} +\frac{1}{16}(\|\nabla
^3d\|^2_{L^2}+\|\nabla^2 u\|^2_{L^2}).\end {array}\eqno (2.6)$$
Similarly,
$$\begin {array} {rl}|I_{1,3}|=&|\langle \partial_x^2u\cdot\nabla d+\partial_xu\cdot\nabla\partial_x d,\partial_x^2 d\rangle|
\\\leq& C\|\triangle d\|^3_{L^3}+\frac{1}{32}\|\nabla^2 u\|^2_{L^2}+C\|\nabla u\|^3_{L^3}\\\leq &C(\|u\|_{L^p}^\frac{2p}{p-3}\|u\|_{L^p}^\frac{p}{p-3}
+\|\nabla d\|_{L^p}^\frac{2p}{p-3}\|\nabla
d\|_{L^p}^\frac{p}{p-3})\\&+\frac{1}{16}(\|\nabla
^3d\|^2_{L^2}+\|\nabla^2 u\|^2_{L^2}),\end {array}\eqno (2.7)$$
$$\begin {array} {rl}|I_{1,4}|&=|\langle \partial_x^1 (|\nabla d|^2d),\partial_x\triangle d\rangle|\\&\leq C(\|\nabla d\|^{6}_{L^6}+\|\nabla^2 d\|^{3}_{L^3})
+\frac{1}{32}\|\nabla\triangle d\|^2_{L^2}\\&\leq C\|\nabla^2
d\|^{3}_{L^3}+\frac{1}{32}\|\nabla\triangle d\|^2_{L^2}\\&\leq
C\|\nabla d\|_{L^p}^\frac{2p}{p-3}\|\nabla
d\|_{L^p}^\frac{p}{p-3}+\frac{1}{16}\|\nabla ^3d\|^2_{L^2}.\end
{array}\eqno (2.8)$$ Taking the above estimates (2.5)-(2.8) in (2.4)
for $|\alpha|=1$, we arrive at
$$\begin {array}{rl}&\frac{d}{dt}\|(\nabla u,\nabla^2 d)\|^2_{L^2}+C\| (\nabla^2 u,\nabla^3 d)\|^2_{L^2}
\\ \leq &C(\|u\|_{L^p}^\frac{2p}{p-3}\|u\|_{L^p}^\frac{p}{p-3}
+\|\nabla d\|_{L^p}^\frac{2p}{p-3}\|\nabla
d\|_{L^p}^\frac{p}{p-3}).\end {array}\eqno (2.9)$$

Next we derive an estimate for $\|u\|^{p}_{L^p}$ and $\|\nabla
d\|^{p}_{L^p}$. First of all, we multiply (1.1a) by $|u|^{p-2}u$ and
integrate over $\mathbb R^3$ to obtain that
$$\begin {array} {rl}&\frac{1}{p}\frac{d}{dt}\|u\|^p_{L^p}+\int_{\mathbb R^3} |u|^{p-2}|\nabla u|^2dx
+\frac{1}{4}(p-2)\int_{\mathbb R^3}|u|^{p-4}|\nabla(|u|^2)|^2dx
\\=&-\langle\nabla P, |u|^{p-2}u\rangle-\langle\triangle d\nabla d, |u|^{p-2}u\rangle.\end {array}\eqno (2.10)$$
Observe that
$$\triangle d\cdot\nabla d=\nabla\cdot(\nabla d\odot\nabla
d-\frac{1}{2}|\nabla d|^2I),\eqno (2.11)$$ where $\nabla
d\odot\nabla d$ denotes the $3\times3$ matrix whose $(i,j)-$the
entry is given by $\partial_id\cdot\partial_jd$ for $1\le i,j\le 3$.
And taking $div$ to (1.1a), we arrive at
$$\triangle P=-div div(u\otimes u+\nabla d\odot \nabla d-\frac{1}{2}|\nabla d|^2I).$$
An application
of the $L^p$-estimate of elliptic systems to the above equation,
there exists $\bar P(t)$ such that
$$\begin {array} {rl}&\int_{\mathbb R^3}|P-\bar P(t)|^{\frac{p+2}{2}}dx\\ \leq& C\int_{\mathbb R^3}|u|^{p+2}+|\nabla d|^{p+2}dx\\
\leq &C\{(\int_{\mathbb R^3}|u|^pdx)^\frac{p-1}{p}(\int_{\mathbb
R^3} |u|^{3p}dx)^\frac{1}{p}\\&+(\int_{\mathbb R^3}|\nabla
d|^pdx)^\frac{p-1}{p}(\int_{\mathbb R^3} |\nabla
d|^{3p}dx)^\frac{1}{p}\}\\\leq& C\{(\int_{\mathbb
R^3}|u|^pdx)^\frac{p-1}{p}(\int_{\mathbb R^3} |u|^{p-2}|\nabla
u|^2dx)^\frac{3}{p}\\&+(\int_{\mathbb R^3}|\nabla
d|^pdx)^\frac{p-1}{p}(\int_{\mathbb R^3} |\nabla d|^{p-2}|\nabla^2
d|^2dx)^\frac{3}{p}\}\\ \leq&C\{(\int_{\mathbb
R^3}|u|^pdx)^\frac{p-1}{p-3}+(\int_{\mathbb R^3}|\nabla
d|^pdx)^\frac{p-1}{p-3}\}\\&+\frac{1}{32}(\int_{\mathbb R^3}
|u|^{p-2}|\nabla u|^2dx+\int_{\mathbb R^3} |\nabla d|^{p-2}|\nabla^2
d|^2dx)\\ \leq&
C\{\|u\|_{L^p}^{\frac{2p}{p-3}}\|u\|_{L^p}^p+\|\nabla
d\|_{L^p}^{\frac{2p}{p-3}}\|\nabla
d\|_{L^p}^p\}\\&+\frac{1}{32}(\int_{\mathbb R^3} |u|^{p-2}|\nabla
u|^2dx+\int_{\mathbb R^3} |\nabla d|^{p-2}|\nabla^2 d|^2dx).\end
{array}$$ Thus we have
$$\begin {array} {rl}&|\langle\nabla P, |u|^{p-2}u\rangle|=|\langle P-\bar P(t), div (|u|^{p-2}u)\rangle|
\\ \leq & C\int_{\mathbb R^3}|\nabla |u|^2||u|^{p-4}|u||P-\bar P(t)|dx\\ \leq& \int_{\mathbb R^3}\frac{1}{32}|\nabla |u|^2|^2|u|^{p-4}+C|u|^{p+2}+|P-\bar P(t)|^\frac{p+2}{2}dx
\\ \leq& C\{\|u\|_{L^p}^{\frac{2p}{p-3}}\|u\|_{L^p}^p+\|\nabla
d\|_{L^p}^{\frac{2p}{p-3}}\|\nabla
d\|_{L^p}^p\}\\&+\frac{1}{32}(\int_{\mathbb R^3}|\nabla
|u|^2|^2|u|^{p-4}+ |u|^{p-2}|\nabla u|^2dx+\int_{\mathbb R^3}
|\nabla d|^{p-2}|\nabla^2 d|^2dx).\end {array}$$ Similarly, applying
(2.11) and by Sobolev's inequality and Cauchy's inequality, we get
$$\begin {array} {rl}&|\langle\triangle d\cdot\nabla d, |u|^{p-2}u\rangle|
=|\langle div (\nabla d\odot\nabla d-\frac{1}{2}|\nabla d|^2I),
|u|^{p-2}u\rangle|
\\ \leq& C\{\|u\|_{L^p}^{\frac{2p}{p-3}}\|u\|_{L^p}^p+\|\nabla
d\|_{L^p}^{\frac{2p}{p-3}}\|\nabla
d\|_{L^p}^p\}\\&+\frac{1}{32}(\int_{\mathbb R^3}|\nabla
|u|^2|^2|u|^{p-4}+ |u|^{p-2}|\nabla u|^2dx+\int_{\mathbb R^3}
|\nabla d|^{p-2}|\nabla^2 d|^2dx).\end {array}$$ Putting the above
two inequalities into (2.10) we have
$$\begin {array} {rl}&\frac{d}{dt}\|u\|^p_{L^p}+\int_{\mathbb R^3} |u|^{p-2}|\nabla u|^2dx
+\frac{1}{2}\int_{\mathbb R^3}|u|^{p-4}|\nabla(|u|^2)|^2dx
\\ \leq& C\{\|u\|_{L^p}^{\frac{2p}{p-3}}\|u\|_{L^p}^p+\|\nabla
d\|_{L^p}^{\frac{2p}{p-3}}\|\nabla
d\|_{L^p}^p\}+\frac{1}{32}\int_{\mathbb R^3} |\nabla
d|^{p-2}|\nabla^2 d|^2dx.\end {array}\eqno (2.12)$$

Next differentiating (1.1b) with respect to $x$, we have
$$\partial_x d_t-\triangle\partial_x d=\partial_x(|\nabla d|^2d-u\cdot\nabla d).\eqno (2.13)$$
We multiply (2.13) by $|\nabla d|^{p-2}\partial_x d$ and integrate
over $\mathbb R^3$ to obtain that
$$\begin {array}{rl}&\frac{1}{p}\frac{d}{dt}\|\nabla d\|^p_{L^p}+\int_{\mathbb R^3}|\nabla d|^{p-2}|\nabla^2
d|^2dx+\frac{1}{4}(p-2)
\int_{\mathbb R^3}|\nabla d|^{p-4}|\partial_x (|\nabla d|^2)|^2dx\\
=& \langle |\nabla d|^{p-2}\partial_x d,\partial_x (|\nabla
d|^2d-u\cdot\nabla
d)\rangle\\
=&\|\nabla d\|^{p+2}_{L^{p+2}} +\langle |\nabla d|^{p-2}\partial_x
d,2\nabla d \nabla\partial_x dd\rangle +\langle \nabla\cdot(|\nabla
d|^{p-2}\partial_x
d),u\cdot\nabla d\rangle\\
\leq& C(\|\nabla
d\|^{p+2}_{L^{p+2}}+\|u\|^{p+2}_{L^{p+2}})+\frac{1}{64}\int_{\mathbb
R^3}|\nabla d|^{p-2}|\nabla^2 d|^2dx\\ \leq&
C\{\|u\|_{L^p}^{\frac{2}{p-3}}\|u\|_{L^p}^p+\|\nabla
d\|_{L^p}^{\frac{2}{p-3}}\|\nabla
d\|_{L^p}^p\}\\&+\frac{1}{32}(\int_{\mathbb R^3} |u|^{p-2}|\nabla
u|^2dx+\int_{\mathbb R^3} |\nabla d|^{p-2}|\nabla^2 d|^2dx).\end
{array}\eqno (2.14)$$ Combining (2.12) and (2.14) gives
$$\begin {array} {rl}&\frac{d}{dt}(\|u\|^p_{L^p}+\|\nabla d\|^p_{L^p})+\int_{\mathbb R^3} |u|^{p-2}|\nabla u|^2+|\nabla d|^{p-2}|\nabla^2 d|^2dx
\\ \leq& C\{\|u\|_{L^p}^{\frac{2p}{p-3}}\|u\|_{L^p}^p+\|\nabla
d\|_{L^p}^{\frac{2p}{p-3}}\|\nabla d\|_{L^p}^p\}.\end {array}\eqno
(2.15)$$

Now by (1.5) we have $q=\frac{2p}{p-3}$, thus
$\|u\|_{L^p}^{\frac{2}{p-3}}$ and $\|\nabla
d\|_{L^p}^{\frac{2}{p-3}}$ belong to $L^1[0,T]$. For $p$, we divide
into
two case:\\

{\bf Case 1.} If $4\leq p<+\infty$, then $\frac{p}{p-3}\leq p$. Thus
combining (2.9) and (2.15) gives
$$\begin {array}{rl}&\frac{d}{dt}(\|(\nabla u,\nabla^2 d)\|^2_{L^2}+\|(u,\nabla d)\|^p_{L^p})+C\| (\nabla^2 u,\nabla^3 d)\|^2_{L^2}
\\ \leq &C(\|u\|_{L^p}^\frac{2p}{p-3}\|u\|_{L^p}^p
+\|\nabla d\|_{L^p}^\frac{2p}{p-3}\|\nabla d\|_{L^p}^p+1).\end
{array}$$ By Gronwall's inequality, we have
$$\begin {array}{rl}&\|(\nabla u,\nabla^2 d)(t,\cdot)\|^2_{L^2}+\|(u,\nabla d)(t,\cdot)\|^p_{L^p}+C\int_0^t\| (\nabla^2 u,\nabla^3
d)(s,\cdot)\|^2_{L^2}ds
\\ \leq &C\exp[C\int_0^t(\|u(s,\cdot)\|_{L^p}^\frac{2p}{p-3}
+\|\nabla d(s,\cdot)\|_{L^p}^\frac{2p}{p-3}+1)ds]<+\infty,\end
{array}\eqno (2.16)$$ for all $t\in[0,T]$.\\

 {\bf Case 2.} If
$3<p<4$. Multiplying (2.15) by $(\|u\|^p_{L^p}+\|\nabla
d\|^p_{L^p})^\frac{3}{p-3}$, we obtain
$$\begin {array} {rl}&\frac{d}{dt}(\|u\|^p_{L^p}+\|\nabla d\|^p_{L^p})^\frac{p}{p-3}\leq C(\|u\|_{L^p}^{\frac{2p}{p-3}}+\|\nabla
d\|_{L^p}^{\frac{2p}{p-3}})(\|u\|^p_{L^p}+\|\nabla
d\|^p_{L^p})^\frac{p}{p-3}.\end {array}\eqno (2.17)$$ We add (2.17)
to (2.9) to obtain
$$\begin {array}{rl}&\frac{d}{dt}[\|(\nabla u,\nabla^2 d)\|^2_{L^2}+(\|u\|^p_{L^p}+\|\nabla d\|^p_{L^p})^\frac{p}{p-3}]+C\| (\nabla^2 u,\nabla^3 d)\|^2_{L^2}
\\ \leq &C(\|u\|_{L^p}^{\frac{2p}{p-3}}+\|\nabla
d\|_{L^p}^{\frac{2p}{p-3}})(\|u\|^p_{L^p}+\|\nabla
d\|^p_{L^p})^\frac{p}{p-3}.\end {array}$$ By Gronwall's inequality,
we get
$$\begin {array}{rl}&\|(\nabla u,\nabla^2 d)(t,\cdot)\|^2_{L^2}+(\|u(t,\cdot)\|^p_{L^p}+\|\nabla d(t,\cdot)\|^p_{L^p})^\frac{p}{p-3}\\&+C\int_0^t\| (\nabla^2 u,\nabla^3
d)(s,\cdot)\|^2_{L^2}ds  \leq C\exp[C\int_0^t\|(u,\nabla
d)(s,\cdot)\|_{L^p}^\frac{2p}{p-3} ds]<+\infty,\end {array}\eqno
(2.18)$$ for all $t\in[0,T]$.

Next for $|\alpha|=3$. For $I_{3,1}$, we need to use the following
Moser-type inequality (see [17, p. 43]):
$$\|D^s(fg)\|_{L^2}\leq C(\|g\|_{L^\infty}\|\nabla^sf\|_{L^2}+\| f\|_{L^\infty}\|\nabla^{s}g\|_{L^2}).\eqno(2.19) $$
Thus we have
$$\begin {array}{rl}|I_{3,1}|&=|\langle \partial_x^2div (u\otimes u),\partial_x^4 u\rangle|
\\&\leq C\|\nabla^3(u\otimes u)\|^2_{L^2}+\frac{1}{16}\|\nabla ^4u\|^2_{L^2}\\&\leq C\|u\|^2_{L^\infty}\|\nabla^3u\|^2_{L^2}+\frac{1}{16}\|\nabla ^4u\|^2_{L^2}
\\&\leq C\{\|\nabla u\|^{\frac{5}{6}}_{L^2}\|\nabla^4 u\|^{\frac{1}{6}}_{L^2}\|\nabla u\|^{\frac{1}{3}}_{L^2}\|\nabla^4
u\|^{\frac{2}{3}}_{L^2}\}^2 +\frac{1}{16}\|\nabla
^4u\|^2_{L^2}\\&\leq C\|\nabla u\|^{14}_{L^2}+\frac{1}{8}\|\nabla
^4u\|^2_{L^2}.\end {array}\eqno (2.20)$$ For $I_{3,2}$, we apply
(2.11) and (2.19) to obtain that
$$\begin {array}{rl}|I_{3,2}|
&=|\langle \partial_x^2 \nabla\cdot(\nabla d\otimes\nabla d-\frac{1}{2}|\nabla d|^2I),\partial_x^4 u\rangle|
\\&\leq C\|\nabla^3(\nabla d\otimes\nabla d-\frac{1}{2}|\nabla d|^2I)\|^2_{L^2}+\frac{1}{8}\|\nabla^4u\|^2_{L^2}
\\&\leq C\|\nabla d\|^2_{L^\infty}\|\nabla^4 d\|^2_{L^2}+\frac{1}{8}\|\nabla ^4u\|^2_{L^2}
\\&\leq C\{\|\nabla^2 d\|^{\frac{5}{6}}_{L^2}\|\nabla^5 d\|^{\frac{1}{6}}_{L^2}\|\nabla^2
d\|^{\frac{1}{3}}_{L^2}\|\nabla^5 d\|^{\frac{2}{3}}_{L^2}\}^2
+\frac{1}{8}\|\nabla ^4u\|^2_{L^2}\\&\leq C\|\nabla^2
d\|^{14}_{L^2}+\frac{1}{8}(\|\nabla ^4u\|^2_{L^2}+\|\nabla
^5d\|^2_{L^2}).\end {array}\eqno(2.21)$$ Similar in the proof of
(2.21), $I_{3,3}, I_{3,4}$ can be bounded as follows:
$$\begin {array}{rl}|I_{3,3}|&=|\langle \partial_x^3 (u\cdot\nabla d),\partial_x^5 d\rangle|\\
 &\leq C\|\partial_x^3 (u\cdot\nabla d)\|^2_{L^2}
+\frac{1}{16}\|\nabla ^5d\|^2_{L^2}\\&\leq C(\|\nabla
d\|^2_{L^\infty}\|\nabla^3 u\|^2_{L^2}+\|u\|^2_{L^\infty}\|\nabla^4
d\|^2_{L^2})+\frac{1}{16}\|\nabla ^5d\|^2_{L^2}
\\&\leq C\{\|\nabla^2 d\|^{\frac{5}{6}}_{L^2}\|\nabla^5 d\|^{\frac{1}{6}}_{L^2}\|\nabla
u\|^{\frac{1}{3}}_{L^2}\|\nabla^4 u\|^{\frac{2}{3}}_{L^2}\\&\ \ \ \
+\|\nabla u\|^{\frac{5}{6}}_{L^2}\|\nabla^4
u\|^{\frac{1}{6}}_{L^2}\|\nabla^2 d\|^{\frac{1}{3}}_{L^2}\|\nabla^5
d\|^{\frac{2}{3}}_{L^2}\}^2 +\frac{1}{16}\|\nabla ^5d\|^2_{L^2}
\\&\leq C(\|\nabla^2
d\|^{14}_{L^2}+\|\nabla u\|^{14}_{L^2})+\frac{1}{8}(\|\nabla
^4u\|^2_{L^2}+\|\nabla ^5d\|^2_{L^2}),\end {array}\eqno(2.22)$$ and
$$\begin {array}{rl}|I_{3,4}|&=|\langle \partial_x^3 (|\nabla d|^2d),\partial_x^5 d\rangle|\\ &\leq C\|\partial_x^3 (|\nabla d|^2d)\|^2_{L^2}
+\frac{1}{32}\|\nabla ^5d\|^2_{L^2}\\&\leq C(\|\nabla^3 |\nabla
d|^2\|^2_{L^2}+\||\nabla d|^2\|^2_{L^\infty}\|\nabla^3
d\|^2_{L^2})+\frac{1}{32}\|\nabla ^5d\|^2_{L^2}\\&\leq C(\|\nabla
d\|^2_{L^\infty}\|\nabla^4 d\|^2_{L^2}+\|\triangle
d\|^2_{L^\infty}\|\nabla^3 d\|^2_{L^2})+\frac{1}{32}\|\nabla
^5d\|^2_{L^2}
\\&\leq C(\|\nabla^2
d\|^{14}_{L^2}+\{\|\nabla^2 d\|^\frac{1}{2}_{L^2}\|\nabla^5
d\|^\frac{1}{2}_{L^2}\|\nabla^2 d\|^\frac{2}{3}_{L^2}\|\nabla^5
d\|^\frac{1}{3}_{L^2}\}^2)+\frac{1}{16}\|\nabla ^5d\|^2_{L^2}\\&\leq
C\|\nabla^2 d\|^{14}_{L^2}+\frac{1}{8}\|\nabla ^5d\|^2_{L^2}.\end
{array}\eqno(2.23)$$ Putting (2.20)-(2.23) into (2.4) for
$|\alpha|=3$ and by (2.16) and (2.18), we arrive at
$$\begin {array} {rl}&\frac{d}{dt}(\|\nabla^3u\|^2_{L^2}+\|\nabla^4 d\|^2_{L^2})+ \|\nabla
^4u\|^2_{L^2}+\|\nabla ^5d\|^2_{L^2}\\&\leq C(\|\nabla^2
d\|^{14}_{L^2}+\|\nabla u\|^{14}_{L^2})\leq +\infty.\end {array}$$
Integrating the above inequality with respect to time from 0 to
$t\in [0,T]$ and by (2.3),
 we conclude $u(T,\cdot)\in H^3(\mathbb R^3)$ and $d(T,\cdot)-a\in
H^4(\mathbb R^3)$. Thus the proof of Theorem 1.1 is completed.\\
\\{\bf Acknowledge}\\

The authors' research was supported Supported by National Natural
Science Foundation of China-NSAF (Grant No. 10976026) and by
National Natural Science Foundation of China-NSAF (Grant No.
11271305).


\begin{thebibliography}{aa}
\footnotesize
 \bibitem{BKM}J. T. Beale, T. Kato, A. Majda, Remarks on the breakdown of smooth
solutions for the 3-D Euler equation. Commun. Math. Phys. 94: 61-66.


\footnotesize \vspace{-0.3cm}
    \bibitem{D}P. G. DE Gennes, The Physics of Liquid Crystals. Oxford,
1974.



\footnotesize \vspace{-0.3cm}
    \bibitem{E}J. L. Ericksen, Hydrostatic theory of liquid crystal. Arch.
Ration. Mech. Anal. 9 (1962), 371-378.

\footnotesize \vspace{-0.3cm}
    \bibitem{ESS}L. Escauriaza, G. Seregin, V. Sver\'ak,
    $L_{3,\infty}$ solutions of the Navier-Stokes equations and
    backward uniqueness. Russ. Math. Surv. 58 (2003), 211-250.


\footnotesize \vspace{-0.3cm}
    \bibitem{H}M. C. Hong, Global existence of solutions of the simplified
Ericksen-Leslie system in $\mathbb R^2$. Calc. Var. Partial
Differential Equations 40 (2011), 15-36.


 \footnotesize \vspace{-0.3cm}
    \bibitem{XH}M. C. Hong, Z. P. Xin, Global existence of solutions of the
Liquid Crystal flow for the Oseen-Frank model in $\mathbb R^2$. Adv.
Math., 231 (2012), 1364-1400.


\footnotesize \vspace{-0.3cm}
    \bibitem{HW}T. Huang, C. Y. Wang, Blow up criterion for nematic liquid
crystal flows.  Comm. Partial Differ. Equ. 37 (2012), 875-884.

\footnotesize \vspace{-0.3cm}
    \bibitem{L1}O. Ladyzhenskaya, Uniqueness and smoothness of
    generalized solutions of Navier-Stokes equations. Zap. Nauv cn.
    Sem. Leningrad. Otdel. Mat. Inst. Steklov. (LOMI) 5 (1967),
    169-185.



 \footnotesize \vspace{-0.3cm}
    \bibitem{L2} J. Leray, Sur le mouvement d'un liquide visqueux
    emplissant l'espace. Acta Math. 63 (1934) 193-248.


\footnotesize \vspace{-0.3cm}
    \bibitem{L3}F. M. Leslie, Some constitutive equations for liquid
crystals. Arch. Ration. Mech. Anal. 28 (1962), 265-283.



\footnotesize \vspace{-0.3cm}
    \bibitem{L4}F. H. Lin, Nonlinear theory of defects in nematic liquid
crystals; phase transition and flow phenomena. Commun. Pure Appl.
Math. 42 (1989), 789-814.




 \footnotesize \vspace{-0.3cm}
    \bibitem{LL1}F. H. Lin, C. Liu, Nonparabolic dissipative systems modeling
the flow of liquid crystals. Commun. Pure Appl. Math. 48 (1995),
501-537.

 \footnotesize \vspace{-0.3cm}
    \bibitem{LL2}F. H. Lin, C. Liu, Partial regularity of the dynamic system
modeling the flow of liquid crystals. Discrete Contin. Dyn. Syst. 2
(1996), 1-22.




  \footnotesize \vspace{-0.3cm}
    \bibitem{LLW1}F. H. Lin, J. Lin, C. Y. Wang, Liquid crystal flows in two
dimensions. Arch. Rational Mech. Anal. 197 (2010), 297-336.

  \footnotesize \vspace{-0.3cm}
    \bibitem{LLW2}F. H. Lin, J. Lin, C. Y. Wang, On the uniqueness of heat flow
of harmonic maps and hydrodynamic flow of nematic liquid crystals.
Chinese Ann. Math. 31B: 921-928.

\footnotesize \vspace{-0.3cm}
    \bibitem{LZ}Q. Liu, J. Zhao, Logarithmical Blow-up Criteria for the Nematic Liquid Crystal
    Flows, printed.



\footnotesize \vspace{-0.3cm}
    \bibitem{M}A. Majda, Compressible fluid flow and system of conservation
laws in several space variables. Applied Mathematical Sciences, 53,
Springer- Verlag, NewYork, 1984.

\footnotesize \vspace{-0.3cm}
    \bibitem{P}G. Prodi, Un teorema di unicit per le equationi di
    Navier-Stokes. Ann. Mat. Pura Appl. (4) 48 (1959), 173-182.

\footnotesize \vspace{-0.3cm}
    \bibitem{S} J. Serrin, The initial value problem for the
    Navier-Stokes equations. Nonlinear Problems (Proc. Sympos., Madison,
    Wis.). University of Wisconsis Press, Madison, 69-98, 1963.


 \footnotesize \vspace{-0.3cm}
    \bibitem{XZ}X. Xu and Z. Zhang, Global regularity and uniqueness of weak
solution for the 2-D liquid crystal flows. J. Differential Equations
252 (2012), 1169-1181.













  \end{thebibliography}
\end{document}